\documentclass[a4paper]{amsart}

\usepackage{longtable}
\usepackage[ruled, lined]{algorithm2e}
\usepackage[dvips, dvipsnames, usenames]{color}

\usepackage{appendix}

\makeatletter

\newcommand{\Aut}{\operatorname{Aut}}

\newcommand{\atlas}{\textsf{ATLAS}}
\newcommand{\gap}{\textsf{GAP}}
\newcommand{\Alt}{\mathbb{A}}
\newcommand{\Sym}{\mathbb{S}}

\newcommand{\oc}{{\mathcal O}}
\newcommand{\M}{{\mathcal M}}

\newcommand\s{\mathbb S}

\providecommand{\\}{\\}

\numberwithin{equation}{section} 
\numberwithin{figure}{section} 
  \theoremstyle{plain}
  \newtheorem{thm}{Theorem}[section]
  \theoremstyle{plain}
  \newtheorem{lem}[thm]{Lemma}
  \theoremstyle{remark}
  \newtheorem{rem}[thm]{Remark}
  \theoremstyle{plain}
  
  \theoremstyle{plain}
  
  \theoremstyle{plain}

\newcounter{maint}

\newtheorem{step-u352}{Step}
\newtheorem{step-l322}{Step}
\newtheorem{step-l332}{Step}
\newtheorem{step-u33}{Step}
\newtheorem{step-2A8}{Step}
\newtheorem{step-A5xD5}{Step}

\newtheorem{step-hs}{Step}
\newtheorem{step-j2}{Step}
\newtheorem{step-co1}{Step}
\newtheorem{step-co2}{Step}
\newtheorem{step-co3}{Step}
\newtheorem{step-mcl}{Step}
\newtheorem{step-suz}{Step}
\newtheorem{step-j1}{Step}
\newtheorem{step-j3}{Step}
\newtheorem{step-ru}{Step}
\newtheorem{step-on}{Step}
\newtheorem{step-ly}{Step}
\newtheorem{step-j4}{Step}
\newtheorem{step-he}{Step}
\newtheorem{step-hn}{Step}
\newtheorem{step-th}{Step}
\newtheorem{step-fi22}{Step}
\newtheorem{step-fi23}{Step}
\newtheorem{step-fi24'}{Step}
\newtheorem{step-b}{Step}
\newtheorem{step-m}{Step}
\newtheorem{step-tits}{Step}

\newcounter{algo}

\theoremstyle{remark}

\newtheorem*{acknowledgement*}{Acknowledgement}
\newtheorem*{acknowledgements*}{Acknowledgements}

\theoremstyle{definition}

\makeatother

\newcommand\toba{{\mathfrak B }}

\def\pf{\begin{proof}}
\def\epf{\end{proof}}

\begin{document}

\title[The \emph{logbook} of Nichols algebras over the sporadic simple 
groups]{The \emph{logbook} of Pointed Hopf algebras over the sporadic simple 
groups}

\author[Andruskiewitsch, Fantino, Gra\~na, Vendramin]{N.  Andruskiewitsch, F.
Fantino, M.  Gra\~na and L. Vendramin}

\thanks{This work was partially supported by
 ANPCyT-Foncyt, CONICET, Ministerio de Ciencia y Tecnolog\'\i a (C\'ordoba)  and Secyt (UNC)}

\address{\noindent N. A., F. F. : Facultad de Matem\'atica, Astronom\'{\i}a y F\'{\i}sica,
Universidad Nacional de C\'ordoba. CIEM -- CONICET. 
Medina Allende s/n (5000) Ciudad Universitaria, C\'ordoba,
Argentina}
\address{\noindent M. G., L. V. : Departamento de Matem\'atica -- FCEyN,
Universidad de Buenos Aires, Pab. I -- Ciudad Universitaria (1428)
Buenos Aires -- Argentina}
\address{\noindent L. V. : Instituto de Ciencias, Universidad de Gral. Sarmiento, J.M. Gutierrez
1150, Los Polvorines (1653), Buenos Aires -- Argentina  }

\address{}

\email{(andrus, fantino)@famaf.unc.edu.ar} \email{(matiasg,
lvendramin)@dm.uba.ar}

\date{\today}

\maketitle

\tableofcontents{}

\section*{Introduction}\label{sect:intro}

In these notes, we give details on the proofs performed with
\textsf{GAP} of the theorems of our paper \cite{AFGV-espo}. Here
we discuss the algorithms implemented for studying Nichols
algebras over non-abelian groups. For each group we refer to files
where the results of our computations are shown. These logs files 
can be downloaded from our webpages

\verb+http://www.mate.uncor.edu/~fantino/afgv-sporadic+

\verb+http://mate.dm.uba.ar/~lvendram/afgv-sporadic+

We believe, and hope, that the details in this work are enough to
guide the reader to repeat and corroborate our calculations.

Throughout the paper, we follow the notations and conventions of
\cite{AFGV-espo}. We write \atlas~ for any of the references
\cite{atlas,AtlasRep1.4,ATLASwww}.




%








\section{Algorithms}\label{sect:algo}
In this section, we explain our algorithms to implement the
techniques presented in \cite[Subsect. 2.2]{AFGV-espo}.

\subsection{Algorithms for type D} \label{subsect:algorithmstipoD}
Algorithm \ref{alg:tipoD} checks if a given conjugacy class $\oc_r^G$
of a finite group $G$ is of type D.

\SetAlgoSkip{bigskip}
\begin{algorithm}[ht]
\dontprintsemicolon \restylealgo{boxed} \caption{Type D}
\label{alg:tipoD} \BlankLine \For{$s\in \oc_r^G$}{
  \If{$(rs)^2 \ne (sr)^2$}{
    Compute the group $H=\langle r,s\rangle$\;
    \If{$\oc_r^H\cap \oc_s^H=\emptyset$}{
    \Return true\tcc*[f]{the class is of type D}
    }
  }
} \Return false\tcc*[f]{the class is not of type D} \BlankLine
\end{algorithm}

\begin{algorithm}[ht]
\dontprintsemicolon \restylealgo{boxed} \caption{Random variation
of Algorithm \ref{alg:tipoD}} \BlankLine \label{alg:tipoD_random}
\ForAll(\tcc*[f]{the number of iterations})
  {$i:1\leq i\leq N$}{
  $x\in G$\tcc*[F]{randomly chosen}
  $s= xrx^{-1}$\;
  \If{$(rs)^2 \ne (sr)^2$}{
    Compute the group $H=\langle r,s\rangle$\;
    \If{$\oc_r^H\cap \oc_s^H=\emptyset$}{
      \Return true\tcc*[f]{the class is of type D}
    }
  }
} \Return false\tcc*[f]{the class is not of type D} \BlankLine
\end{algorithm}

Notice that in Algorithm \ref{alg:tipoD} we need to run over all
the conjugacy class $\oc_r^G$ to look for the element $s$ such
that the conditions of \cite[Prop. 3.4]{AFGV} are satisfied. This
is not always an easy task.  To avoid this problem, we have {\em
the random variation} of Algorithm \ref{alg:tipoD}, that is our
Algorithm \ref{alg:tipoD_random}. The key is to pick randomly an
element $x$ in the group $G$ and to check if $r$ and $s=xrx^{-1}$
satisfy the conditions of \cite[Prop. 3.4]{AFGV}; if not we repeat
the process $N$ times, where $N$ is fixed. This naive variation of
the Algorithm \ref{alg:tipoD} turns out to be very powerful and
allows us to study big sporadic groups such as the Janko group
$J_4$ or the Fischer group $Fi_{24}'$.

For large groups, it is more economical to implement Algorithms
\ref{alg:tipoD} and \ref{alg:tipoD_random} in a recursive way.
Let $G$ be a finite group represented faithfully, for example, as
a permutation group or inside a matrix group over a finite field.
We compute $\oc_{g_1}^G,\dots,\oc_{g_n}^G$, the set of conjugacy
classes of $G$. To decide if these conjugacy classes are of type
D, we restrict the computations to be done inside a nice subgroup
of $G$.

Assume that the list of all maximal subgroups of $G$, up to
conjugacy, is known, namely $\M_1,\M_2,\dots,\M_k$, with
non-decreasing order. Also assume that it is possible to restrict our
good representation of $G$ to every maximal subgroup $\M_i$. Here
we say that a representation is \emph{good} if it allows us to
perform our computations in a reasonable time.


Fix $i \in  \{1, \dots, k\}$. Let $h\in M_i$ and let
$\oc^{\M_i}_h$ be the conjugacy class of $h$ in $\M_i$. Since
$\M_i$ is a subgroup of $G$, the element $h$ belongs to a
conjugacy class of $G$, say $\oc_h^G$. So, if the class
$\oc_h^{\M_i}$ is of type D, then the class $\oc_h^G$ is of type D
too.

\begin{algorithm}[ht]
\BlankLine \dontprintsemicolon \restylealgo{boxed} \caption{Type
D: Using maximal subgroups} \label{alg:tipoD_maximales}
$\oc_{g_1}^G,\dots,\oc_{g_n}^G$ is the set of conjugacy classes of $G$\;
$S=\{1,2,\dots,n\}$\; \BlankLine \ForEach{maximal
subgroup $\M$}{
  Compute $\oc_{h_1}^\M,\dots,\oc_{h_m}^\M$, the set of conjugacy classes of $\M$\;
  \ForEach{$i:1\leq i\leq m$}{
  Identify $\oc_{h_i}^\M$ with a conjugacy class in $G$: $h_i\in \oc_{g_{\sigma(i)}}^G$\;
    \If{$\sigma(i)\in S$}{
      \If{$\oc_{h_i}^\M$ is of type D}{
        Remove $\sigma(i)$ from $S$\;
    \If{$S=\emptyset$}{
          \Return true\tcc*[f]{the group is of type D}
    }
      }
    }
  }
} \ForEach{$j\in S$}{
  \If{$\oc_{g_j}^G$ is of type D}{
    Remove $s$ from $S$\;
    \If{$S=\emptyset$}{
      \Return true\tcc*[f]{the group is of type D}
    }
  }
} \Return $S$\tcc*[f]{conjugacy classes not of type D} \BlankLine
\end{algorithm}

Notice that to implement Algorithm \ref{alg:tipoD_maximales} we
need to have not only a good representation for the group $G$, but
we need to know how to restrict the good representation of $G$ to
all its maximal subgroups. This information appears in the
\atlas~for many of the sporadic simple groups. So, using the
\gap~interface to the \atlas~
we could implement Algorithm \ref{alg:tipoD_maximales} for the sporadic
simple groups.

\subsection{Structure constants}

Some conjugacy classes of involutions are studied 
with \cite[Prop. 1.7, Equation (1.4)]{AFGV-espo}. For example, 
in the proof of Theorem \ref{thm:PSL(5,2)} we claim 
that the conjugacy class 2A of $L_5(2)$ gives only infinite-dimensional 
Nichols algebras. This follows from \cite[Prop. 1.8]{AFGV-espo} 
because $S(\textup{2A},\textup{3A},\textup{3A})=42$. 
\begin{verbatim}
     gap> ct := CharacterTable("L5(2)");;
     gap> ClassNames(ct);;
     gap> ClassMultiplicationCoefficient(ct, ct.2a, ct.3a, ct.3a);
     42
\end{verbatim}

\subsection{Bases for permutation groups}
Let $G$ be a group acting on a set $X$. A subset $B$ of $X$ is
called a \textit{base} for $G$ if the identity is the only element of $G$
which fixes every element in $B$, see \cite[Subsection 3.3]{dm}. In other words, 
\[
\{g\in G\mid g\cdot b=b,\,\text{for all }b\in B\}=1.
\]

\begin{lem}
Let $G$ be a group acting on a set $X$. Let $B$ be a subset of $X$.
The following are equivalent:
\begin{enumerate}
\item $B$ is a base for $G$.
\item For all $g,h\in G$ we have: $g\cdot b=h\cdot b$ for all $b\in B$
implies $g=h$.
\end{enumerate}
\end{lem}

\begin{proof}
If $B$ is a base, then $g\cdot b=h\cdot b\Rightarrow (h^{-1}g)\cdot b=b\Rightarrow h^{-1}g=1\Rightarrow h=g$.
The converse is trivial.
\end{proof}

Let $G$ be a permutation group. 
With the \textsf{GAP} function \texttt{BaseOfGroup} we compute a base for $G$. We use
\texttt{OnTuples} to encode a permutation and \texttt{RepresentativeAction}
to decode the information. This enables us to reduce the size of our log files.

\subsection{An algorithm for involutions} \label{subsect:algo:involution}



We describe here an algorithm used to discard some conjugacy
classes of involutions. In this work, we use this algorithm for
the classes called \textup{2A} in $O_7(3)$, $S_6(2)$, $S_8(2)$,
$Co_2$, $Fi_{22}$, $Fi_{23}$ and $B$.

\smallbreak

Let $\oc$ be one of the classes \textup{2A} in $O_7(3)$, $S_6(2)$,
$S_8(2)$ or $Co_2$, and $g\in \oc$. It is
enough to consider the irreducible representations $\rho$ of the
corresponding centralizer such that $\rho(g)=-1$, 
see \cite[(1.3)]{AFGV-espo}. For these
conjugacy classes, the remaining representations $\rho$ satisfy
$\deg\rho>4$, as can be seen from the character tables. 
By \cite[Lemma 1.3]{AFGV-espo}, we are reduced to
find an involution $x$ such that $gh=hg$, for $h=xgx^{-1}$, and to
compute the multiplicities of the eigenvalues of $\rho(h)$. For
this last task, we use \cite[Remark 1.4]{AFGV-espo}. See the proofs
of Theorems \ref{O7(3)_tipoB}, \ref{thm:S6(2)} and Lemma
\ref{lem:S8(2)_tipoB}, for the classes \textup{2A} of $O_7(3)$,
\textup{2A} of $S_6(2)$ and \textup{2A} of $S_8(2)$, respectively,
and the file \texttt{Co2/2A.log} for \textup{2A} of $Co_2$.


\smallbreak

Let $\oc$ be the class \textup{2A} of $Fi_{22}$; it has $3510$
elements. Let $\rho$ be an irreducible representation of the
corresponding centralizer. The group $O_7(3)$ is a maximal
subgroup of $Fi_{22}$ and the class \textup{2A} of $O_7(3)$ is
contained in the class \textup{2A} of $Fi_{22}$ -- see the file
\texttt{Fi22/2A.log}. By \cite[Lemma 1.5]{AFGV-espo},
$\dim\toba(\oc,\rho)=\infty$.

\smallbreak

Let $\oc$ be the class \textup{2A} of $Fi_{23}$. Let $\rho$ be an
irreducible representation of the corresponding centralizer. The
group $S_8(2)$ is a maximal subgroup of $Fi_{23}$ and the class
\textup{2A} of $S_8(2)$ is contained in the class \textup{2A} of
$Fi_{23}$. By 
\cite[Lemma 1.5]{AFGV-espo}, $\dim\toba(\oc,\rho)=\infty$.

\smallbreak

Let $\oc$ be the class \textup{2A} of $B$. Let $\rho$ be an
irreducible representation of the corresponding centralizer. The
group $Fi_{23}$ is a maximal subgroup of $B$ and the class
\textup{2A} of $Fi_{23}$ is contained in the class \textup{2A} of
$B$ -- see the file \texttt{B/step2.log}. By \cite[Lemma
1.5]{AFGV-espo}, $\dim\toba(\oc,\rho)=\infty$.


\section{Some auxiliary groups}\label{sect:groupsLie}

In this section we study some groups that appear as subquotients
of some sporadic simple groups.

\subsection{The groups $\mathbb A_9$, $\mathbb A_{11}$, $\mathbb A_{12}$, $\mathbb
S_{12}$}\label{subsect:altsim} In \cite{AFGV}, we classify the
conjugacy classes of type D in alternating and symmetric groups.
In Table \ref{tab:altsim}, we list all non-trivial permutations
such that their conjugacy classes are not of type D for the groups
$\mathbb A_9$, $\mathbb A_{11}$, $\mathbb A_{12}$, $\mathbb
S_{12}$.

\begin{table}[ht]
\caption{Classes not of type D in some alternating and
symmetric groups.} \label{tab:altsim}
\begin{center}
\begin{tabular}{|c|c|c|}
\hline Group & Not of type D & Log file\tabularnewline \hline
$\mathbb{A}_{9}$ & $(1\,2\,3)$ &
$\texttt{A9/A9.log}$\tabularnewline \hline $\mathbb{A}_{11}$ &
$(1\,2\,3)$ & $\texttt{A11/A11.log}$\tabularnewline
 & $(1\,2\,3\,4\,5\,6\,7\,8\,9\,10\,11)$ & \tabularnewline
 & $(1\,2\,3\,4\,5\,6\,7\,8\,9\,11\,10)$ & \tabularnewline
\hline $\mathbb{A}_{12}$ & $(1\,2\,3)$ &
$\texttt{A12/A12.log}$\tabularnewline
 & $(1\,2\,3\,4\,5\,6\,7\,8\,9\,10\,11)$ & \tabularnewline
 & $(1\,2\,3\,4\,5\,6\,7\,8\,9\,10\,12)$ & \tabularnewline
\hline $\mathbb{S}_{12}$ & $(1\;2)$ &
$\texttt{S12/S12.log}$\tabularnewline
 & $(1\,2\,3)$ & \tabularnewline
\hline
\end{tabular}
\end{center}
\end{table}

\subsection{The group $L_5(2)$}\label{subsect:L_5(2)}

This group has order $9\,999\,360$. It has 27 conjugacy classes. To study
Nichols algebras over this group we use the representation inside $\Sym_{31}$
given in the \atlas.

\begin{lem}
\label{lem:PSL(5,2)} Every non-trivial conjugacy class of
$L_5(2)$, except \textup{2A} and those with representatives of
order $31$, is of type D.
\end{lem}

\begin{proof}
We perform Algorithm \ref{alg:tipoD}, see the file
\verb+L5(2)/L5(2).log+ for details.
\end{proof}

\begin{thm}
\label{thm:PSL(5,2)} 
The group $L_5(2)$ collapses.
\end{thm}

\begin{proof}
Let $\oc$ be the class \textup{2A}; then
$\dim\mathfrak{B}(\mathcal{O}, \rho)=\infty$ for any irreducible
representation $\rho$ of the corresponding centralizer since
$S(\textup{2A},\textup{3A},\textup{3A})=42$ and \cite[Prop.
1.8]{AFGV-espo} applies. 
Now the
result follows from Lemma \ref{lem:PSL(5,2)} and from \cite[Lemma
1.2]{AFGV-espo} for the conjugacy classes with representatives of
order 31 since these are quasi-real of type $j=2$ and $g^{j^2}\neq
g$.
\end{proof}

\subsection{The group $O_7(3)$}\label{subsect:O_7(3)}

This group has order $4\,585\,351\,680$. It has 58 conjugacy classes. For
computations we use a representation inside $\Sym_{351}$ given in the \atlas.

\begin{lem}
\label{lem:O7(3)_tipoD} Every non-trivial conjugacy class of $O_7(3)$, except \textup{2A},
is of type D.
\end{lem}

\begin{proof}
We perform Algorithm \ref{alg:tipoD_random}, see the file
\verb+O7(3)/O7(3).log+ for details.
\end{proof}

\begin{thm}
\label{O7(3)_tipoB}
The orthogonal group $O_7(3)$ collapses.
\end{thm}

\begin{proof}
By Lemma \ref{lem:O7(3)_tipoD} it remains to study the conjugacy
class \textup{2A}.  For this conjugacy class we use \cite[Lemma 1.3]{AFGV-espo}. 
See the file \verb+O7(3)/2A.log+ for details.
\end{proof}

\subsection{The group $O^+_8(2)$}\label{subsect:O+8(2)}

This group has order $174\,182\,400$. It has 53 conjugacy classes.
For the computations we construct a permutation representation.

\begin{lem}
\label{lem:O8+(2)_tipoD}
Every non-trivial conjugacy class with representative of order distinct from $2$, $3$ is of type D.
\end{lem}

\begin{proof}
We perform Algorithm \ref{alg:tipoD_random}, see the file
\texttt{O8+(2)/O8+2(2).log} for details.
\end{proof}

\begin{thm}
The group $O^+_8(2)$ collapses.
\end{thm}

\begin{proof}
By Lemma \ref{lem:O8+(2)_tipoD} it remains to consider the
conjugacy clases with representative of order 2 or 3. For the five
conjugacy classes of involutions in $O^+_8(2)$ we use \cite[Prop.
1.8]{AFGV-espo}. See Table \ref{tab:O8+(2)_involuciones} for
details.

\begin{table}[ht]
\caption{Involutions in $O^+_8(2)$.}
\label{tab:O8+(2)_involuciones}
\begin{center}
\begin{tabular}{|c|c|c|}
\hline
Class & Size & \\
\hline
\textup{2A} & $1575$ 
& $S(\textup{2A},\textup{3E},\textup{3E})=81$\\
\textup{2B} & $3780$ 
& $S(\textup{2B},\textup{3E},\textup{3E})=108$\\
\textup{2C} & $3780$ 
&  $S(\textup{2C},\textup{3E},\textup{3E})=108$\\
\textup{2D} & $3780$ 
&  $S(\textup{2D},\textup{3E},\textup{3E})=108$\\
\textup{2E} & $56700$ 
&  $S(\textup{2E},\textup{3E},\textup{3E})=486$\\
\hline
\end{tabular}
\end{center}
\end{table}

On the other hand, the conjugacy classes \textup{3A, 3B, 3C, 3D, 3E} of $O^+_8(2)$ are real,
so \cite[Lemma 2.2]{AZ},  cf. \cite[Lemma 1.2]{AFGV-espo}, applies, and the result follows.
\end{proof}

\subsection{The group $O^-_{10}(2)$}\label{subsect:O^-_10(2)}

This group has order $25\,015\,379\,558\,400$. It has 115 conjugacy classes.
For the computations we use the representation inside $\Sym_{495}$ given in the
\atlas.

\begin{lem}
\label{lem:O10-(2)_tipoD} Every non-trivial conjugacy class, except the
conjugacy classes \textup{2A, 3A, 11A, 11B, 33A, 33B, 33C,
33D}, is of type D. 
\end{lem}

\begin{proof}
We perform Algorithm \ref{alg:tipoD_random}, see the file
\verb+O10-(2)/O10-(2).log+ for details.
\end{proof}

\begin{thm}
The group $O^-_{10}(2)$ collapses.
\end{thm}

\begin{proof}
By Lemma \ref{lem:O10-(2)_tipoD} it remains to study the conjugacy
classes \textup{2A, 3A, 11A, 11B, 33A, 33B, 33C, 33D}. Let $\oc$
be the class \textup{2A}; then $\dim\mathfrak{B}(\mathcal{O},
\rho)=\infty$ for any irreducible representation $\rho$ of the
corresponding centralizer since
$S(\textup{2A},\textup{3F},\textup{3F})=243$ and \cite[Prop.
1.8]{AFGV-espo} applies. 
For the class 3A use \cite[Lemma 2.2]{AZ}, since it is a real
conjugacy class.  And for the classes \textup{11A, 11B, 33A, 33B,
33C, 33D} use \cite[Lemma 1.2]{AFGV-espo}, since the classes
\textup{11A, 11B} (resp. \textup{33A, 33B, 33C, 33D}) are
quasi-real of type $j=3$ (resp. $j=4$) with $g^{j^2}\neq g$.
\end{proof}




\subsection{The exceptional group $G_2(4)$}\label{subsect:G2(4)}

In this section we prove that the group $G_2(4)$ collapses.
This group has order $251\,596\,800$. It has 32 conjugacy classes.
In particular, the conjugacy classes with representatives of order
2 or 3 are the following:

\begin{center}
\begin{tabular}{c|c}
Name & Centralizer size\\
\hline
2A & 61440 \\
2B & 3840\\
3A & 60480\\
3B & 180 \\
\end{tabular}
\end{center}

\begin{lem}\label{lem:G2(4)_tipoD}
The conjugacy classes of $G_2(4)$ with representatives of order distinct from $2$, $3$ are of
type D.
\end{lem}

\begin{proof}
We perform Algorithm \ref{alg:tipoD_maximales}, see the file
\verb+G2(4)/G2(4).log+ for details.
\end{proof}

\begin{thm}
\label{thm:G2(4)}
The group $G_2(4)$ collapses.\qed
\end{thm}

\begin{proof}
By Lemma \ref{lem:G2(4)_tipoD}, it remains to study the conjugacy
classes with representatives of order 2 or 3. For the two
conjugacy classes of involutions use \cite[Prop. 1.8]{AFGV-espo},
because $S(\textup{2A},\textup{3B},\textup{3B})=171$ and
$S(\textup{2B}, \textup{3A}, \textup{3A})=126$. 
For the conjugacy
classes with representatives of order 3 use \cite[Lemma
1.2]{AFGV-espo}, because these conjugacy classes are real.
\end{proof}

\subsection{The exceptional groups $G_2(3)$ and
$G_2(5)$}\label{subsect:G2(3):G2(5)}

For the orders and number of conjugacy classes of the groups $G_2(3)$ and
$G_2(5)$ see Table \ref{tab:ChevalleyD}.  We have the following result.

\begin{thm}
\label{thm:G2(3)}
\label{thm:G2(5)}
The groups $G_2(3)$ and $G_2(5)$ are of type D. Hence, they collapse.
\end{thm}

\begin{proof}
We perform Algorithm \ref{alg:tipoD_maximales}, see Table
\ref{tab:ChevalleyD} for the log files.
\end{proof}

\begin{table}[ht]
\caption{Some Chevalley groups of type D.}
\label{tab:ChevalleyD}
\begin{center}
\begin{tabular}{|c|c|c|c|}
\hline
Group & Order & Conjugacy classes & Log file\\
\hline
$G_2(3)$ & $4\,245\,696$ & 23 & \verb+G2(3)/G2(3).log+ \\
\hline
$G_2(5)$ & $5\,859\,000\,000$ & 44 & \verb+G2(5)/G2(5).log+ \\
\hline
\end{tabular}
\end{center}
\end{table}

\subsection{The symplectic group $S_6(2)$}\label{subsect:S6(2)}

This group has order $1\,451\,520$. It has 30 conjugacy classes. For the
computations we use the representation of $S_6(2)$ inside $\Sym_{28}$ given in the
\atlas.

\begin{lem}
\label{lem:S6(2)_tipoD} Every non-trivial conjugacy class of $S_6(2)$, with
the possible exception of \textup{2A, 2B, 3A}, is of type D.
\end{lem}

\begin{proof}
We perform Algorithm \ref{alg:tipoD_random}, see the file
\verb+S6(2)/S6(2).log+ for details.
\end{proof}


\begin{thm}
\label{thm:S6(2)}
The group $S_6(2)$ collapses.
\end{thm}

\begin{proof}
By Lemma \ref{lem:S6(2)_tipoD} it remains to study the classes 2A,
2B, 3A.  For the conjugacy class \textup{2A} use \cite[Lemma
1.3]{AFGV-espo}, see the file \verb+S6(2)/2A.log+. For the conjugacy
class \textup{2B} use \cite[Prop. 1.8]{AFGV-espo}, since
$S(\textup{2B}, \textup{3C}, \textup{3C})=27$. 
For the conjugacy class \textup{3A} use
\cite[Lemma 2.2]{AZ}, since this conjugacy class is real.
\end{proof}

\subsection{The symplectic group $S_8(2)$}\label{subsect:S8(2)}

This group has order $47\,377\,612\,800$. It has 81 conjugacy classes. For the
computations we use the representation of $S_8(2)$ inside $\Sym_{120}$ given in
the \atlas.

\begin{lem}
\label{lem:S8(2)_tipoD} Every non-trivial conjugacy class, except \textup{2A,
2B, 3A}, is of type D.
\end{lem}

\begin{proof}
We perform Algorithm \ref{alg:tipoD_random}, see the file
\verb+S8(2)/S8(2).log+ for details.
\end{proof}


\begin{lem}
\label{lem:S8(2)_tipoB} Let $\oc$ be one of the classes
\textup{2A, 2B, 3A}. Then $\dim\mathfrak{B}(\mathcal{O},
\rho)=\infty$ for any irreducible representation $\rho$ of the
corresponding centralizer.
\end{lem}

\begin{proof}
For the real conjugacy class \textup{3A} use \cite[Lemma 1.2]{AFGV-espo}. For
the conjugacy class \textup{2A} use \cite[Lemma 1.3]{AFGV-espo}, see the file
\verb+S8(2)/2A.log+ for details.  For the conjugacy class \textup{2B} use
\cite[Prop.  1.8]{AFGV-espo}, since
$S(\textup{2B},\textup{3C},\textup{3C})=135$. 
\end{proof}

\subsection{The automorphism group of the Tits group}

In this subsection we study Nichols algebras over the group
$\Aut(^2F_4(2)')\simeq\Aut(^2F_4(2))\simeq{}^2F_4(2)$, the
automorphism group of the Tits group, see \cite{GL}. This group
has order $35\,942\,400$. It has 29 conjugacy classes.

\begin{lem}
\label{lem:2F4(2)_tipoD}
Every non-trivial conjugacy class of $^2F_4(2)$, with the exception of the class {\em 2A}, of size 1755,
is of type D.
\end{lem}

\begin{proof}
We perform Algorithm \ref{alg:tipoD}. See the file
\texttt{2F4(2)/2F4(2).log} for details.
\end{proof}

\begin{thm}
The group $^2F_4(2)$ collapses.
\end{thm}

\begin{proof}
By Lemma \ref{lem:2F4(2)_tipoD} it remains to study the conjugacy
class 2A. For this conjugacy class use \cite[Prop.
1.8]{AFGV-espo}, because
$S(\textup{2A},\textup{3A},\textup{3A})=27$. 
\end{proof}

\subsection{Direct products}\label{subsect:othergroups}

In this subsection we study some direct products of groups that appear as
subgroups or subquotients of the sporadic simple groups.


\subsubsection*{The group $\mathbb{A}_{4}\times G_{2}(4)$}\label{exa:A4xG2(4)}
In this group every non-trivial conjugacy class with representative of order distinct from $2, 3, 6$
is of type D. This follows from Lemma \ref{lem:G2(4)_tipoD} and \cite[Lemma 2.8]{AFGV-espo}.

\subsubsection*{The group $3\times G_{2}(3)$}\label{exa:3xG2(3)} In this group every non-trivial
conjugacy class with representative of order distinct from $3$ is of type D
since every non-trivial conjugacy class of $G_{2}(3)$ is of type D. This follows from Theorem \ref{thm:G2(3)} and \cite[Lemma 2.8]{AFGV-espo}.


\subsubsection*{The group $\mathbb{A}_{9}\times\mathbb{S}_{3}$}\label{exa:A9xS3}
In this group every non-trivial conjugacy class with
representative of order distinct from $2,3$ is of type D. 
See the file \texttt{A9xS3/A9xS3.log} for the computations.




\subsubsection*{The group $\mathbb{S}_5\times L_3(2)$}\label{exa:S5xL3(2)}
In this group every non-trivial conjugacy class with
representative of order distinct from $2, 3, 4, 6, 7$ is of type D.
This was proved with Algorithm \ref{alg:tipoD} -- see the file
\texttt{S5xL3(2)/S5xL3(2).log} for details.

\subsubsection*{The group $\Alt_6\times U_3(3)$}
\label{exa:A6xU3(3)}
In this group every conjugacy class with representative of order 28 or 35 is of
type D. This was proved with Algorithm \ref{alg:tipoD} -- see the file
\texttt{A6xU3(3)/A6xU3(3).log} for details.

\subsubsection*{The group $\mathbb S_5\times \mathbb S_9$}
\label{exa:S5xS9}

In this group every non-trivial conjugacy class with
representative of order distinct from $2, 3, 6$ is of type D. This
was proved with Algorithm 2 -- see the file
\texttt{S5xS9/S5xS9.log} for details.

\section{Proof of Theorem II}\label{sect:largegroups}

In Table \ref{tab:alg2} we list the conjugacy classes of the
sporadic simple groups that are not necessarily of type D. We use
the phrase ``all collapse'' to indicate those groups where all
non-trivial conjugacy classes are of type D. Also we give a
reference about either the algorithms or else the subsections
where these groups are treated.

In the next subsections, we deal with some large sporadic groups.
In order to study these groups we study their maximal subgroups
with Algorithms \ref{alg:tipoD} or \ref{alg:tipoD_random}. The
fusion of the conjugacy classes of the maximal subgroups of a
sporadic group is stored in the \atlas, up to the Monster group
$M$. For this group the fusion of the conjugacy classes is known
only for some of its maximal subgroups; furthermore, the list of
all maximal subgroups of $M$ is not known. In the case of the
Baby Monster group $B$ all the maximal subgroups and the fusions
of conjugacy classes are known except the fusion of the conjugacy
classes of the sixth maximal subgroup.

In each case we split the proof into several steps according to
the corresponding maximal subgroup. We collect in tables and logs
the relevant information.




\begin{table}[ht]
\caption{Conjugacy classes not known of type D.} \label{tab:alg2}
\begin{center}
\begin{tabular}{|l|c|c|}
\hline $G$ & {\bf Conjugacy classes not necessarily of type D}  & {\bf Reference}
\\ \hline $M_{11}$ &  \textup{8A, 8B, 11A, 11B} & Algorithm \ref{alg:tipoD_maximales}
\\ \hline $M_{12}$ &  \textup{11A, 11B} & Algorithm \ref{alg:tipoD_maximales}
\\ \hline $M_{22}$ &  \textup{11A, 11B} & Algorithm \ref{alg:tipoD_maximales}
\\ \hline $M_{23}$ &  \textup{23A, 23B} & Algorithm \ref{alg:tipoD_maximales}
\\ \hline $M_{24}$ &  \textup{23A, 23B} & Algorithm \ref{alg:tipoD_maximales}
\\ \hline $J_{1}$ &  \textup{15A, 15B, 19A, 19B, 19C} & Algorithm \ref{alg:tipoD_maximales}
\\ \hline $J_{2}$ &  \textup{2A, 3A} & Algorithm \ref{alg:tipoD_maximales}
\\ \hline $J_{3}$ &  \textup{5A, 5B, 19A, 19B} & Algorithm \ref{alg:tipoD_maximales}
\\ \hline $Suz$ &  \textup{3A} & Algorithm \ref{alg:tipoD_maximales}
\\ \hline $Ru$ &  \textup{29A, 29B} & Algorithm \ref{alg:tipoD_maximales}
\\ \hline $HS$ &  \textup{11A, 11B} & Algorithm \ref{alg:tipoD_maximales}
\\ \hline $He$ &  \textup{all collapse} & Algorithm \ref{alg:tipoD_maximales}
\\ \hline $McL$ &  \textup{11A, 11B} &Algorithm \ref{alg:tipoD_maximales}
\\ \hline $Co_{3}$ &  \textup{23A, 23B} & Algorithm \ref{alg:tipoD_maximales}
\\ \hline $Co_{2}$ &  \textup{2A, 23A, 23B} & Algorithm \ref{alg:tipoD_maximales}
\\ \hline $O'N$ &  \textup{31A, 31B}  & Algorithm \ref{alg:tipoD_maximales}
\\ \hline $Fi_{22}$   & \textup{2A, 22A, 22B} & Algorithm \ref{alg:tipoD_maximales}
\\ \hline $T= {}^2F_4(2)'$ &  \textup{2A} & Algorithm \ref{alg:tipoD_maximales}
\\ \hline $Co_1$ &  \textup{3A, 23A, 23B} & \S\ref{subsect:Co1}
\\ \hline $Fi_{23}$ &  \textup{2A, 23A, 23B} & \S\ref{subsect:Fi23}
\\ \hline $HN$ &  \textup{all collapse} & \S\ref{subsect:HN}
\\ \hline $Th$ &  \textup{all collapse} & \S\ref{subsect:Th}
\\ \hline $Ly$ &  \textup{33A, 33B, 37A, 37B, 67A, 67B, 67C} & \S\ref{subsect:Ly}
\\ \hline $J_4$ &  \textup{29A, 37A, 37B, 37C, 43A, 43B, 43C} & \S\ref{subsect:J4}
\\ \hline $Fi'_{24}$    & \textup{23A, 23B, 27B, 27C} & \S\ref{subsect:Fi24'}
\\ \hline       & \textup{29A, 29B, 33A, 33B, 39C, 39D} &
\\ \hline $B$   & \textup{2A, 16C, 16D, 32A} & \S\ref{subsect:B}
\\ \hline   & \textup{32B, 32C, 32D, 34A} &
\\ \hline   & \textup{46A, 46B, 47A, 47B} &
\\ \hline $M$   & \textup{32A, 32B, 41A, 46A, 46B} & \S\ref{subsect:M}
\\ \hline   & \textup{47A, 47B, 59A} &
\\ \hline   & \textup{59B, 69A, 69B, 71A, 71B} &
\\ \hline   & \textup{87A, 87B, 92A, 92B, 94A, 94B} &

\\ \hline
\end{tabular}
\end{center}
\end{table}

\begin{table}[ht]
\caption{Log files for the sporadic  groups studied with Algorithm \ref{alg:tipoD_maximales}}
\label{tab:alg2_logs}
\begin{center}
\begin{tabular}{|l|c||l|c|}
\hline $G$ & {\bf Log file} &  $G$ & {\bf Log file}
\\ \hline $M_{11}$ & \texttt{M11/M11.log} & $Ru$ & \texttt{Ru/Ru.log}
\\ \hline $M_{12}$ & \texttt{M12/M12.log} & $HS$ & \texttt{HS/HS.log}
\\ \hline $M_{22}$ & \texttt{M22/M22.log}& $He$ & \texttt{He/He.log}
\\ \hline $M_{23}$ & \texttt{M23/M23.log}& $McL$ & \texttt{McL/McL.log}
\\ \hline $M_{24}$ & \texttt{M24/M24.log}& $Co_{3}$ & \texttt{Co3/Co3.log}
\\ \hline $J_{1}$ & \texttt{J1/J1.log}& $Co_{2}$ &\texttt{Co2/Co2.log}
\\ \hline $J_{2}$ & \texttt{J2/J2.log}& $O'N$ & \texttt{ON/ON.log}
\\ \hline $J_{3}$ & \texttt{J3/J3.log}& $Fi_{22}$ & \texttt{Fi22/Fi22.log}
\\ \hline $Suz$ & \texttt{Suz/Suz.log}& $T$ & \texttt{T/T.log}
\\ \hline
\end{tabular}
\end{center}
\end{table}

\subsection{The Lyons group $Ly$}\label{subsect:Ly}

\begin{step-ly}
The maximal subgroup $\M_1 \simeq G_{2}(5)$.
\end{step-ly}

In $G_2(5)$ every non-trivial conjugacy class is of type D, see Theorem \ref{thm:G2(5)}.
Therefore, the conjugacy
classes \textup{2A, 3A, 3B, 4A, 5A, 5B, 6A, 6B, 6C, 7A, 8A, 8B,
10A, 10B, 12A, 12B, 15A, 15B, 15C, 20A, 21A, 21B, 24A, 24B, 24C,
25A, 30A, 30B, 31A, 31B, 31C, 31D, 31E} of $Ly$ are of type D.

\begin{step-ly}
The maximal subgroup $\M_4 \simeq 2.\mathbb{A}_{11}$.
\end{step-ly}

 Consider the short exact sequence
$1\to2\to2.\mathbb{A}_{11}\to\mathbb{A}_{11}\to 1$. From Table
\ref{tab:altsim}, every non-trivial conjugacy class of
$\mathbb{A}_{11}$ is of type D except the class of the 3-cycles
and the classes of the 11-cycles. Thus, every non-trivial
conjugacy class in $2.\mathbb{A}_{11}$ with representative of
order distinct from $2$, $3$, $6$, $11$, $22$ is of type D, by
\cite[Lemma 2.7]{AFGV-espo}. Hence the conjugacy classes
\textup{9A, 14A, 18A, 28A, 40A, 40B, 42A, 42B} of $Ly$ are of type
D.

\begin{step-ly}
The maximal subgroup $\M_6\simeq 3^5:(2\times M_{11})$.
\end{step-ly}

 We construct a permutation representation of $\M_6$ and
apply Algorithm \ref{alg:tipoD} in this maximal subgroup. We check
that all conjugacy classes of $\M_6$ with representative of order
$11, 22$ are of type D -- see the file \texttt{Ly/step3.log}. By
the fusion of the conjugacy classes, the conjugacy classes
\textup{11A, 11B, 22A, 22B} of $Ly$ are of type D.

\begin{rem}
Not necessarily of type D: \textup{33A, 33B, 37A, 37B, 67A, 67B,
67C}.
\end{rem}

\subsection{The Thompson group $Th$}\label{subsect:Th}

See the file \texttt{Th/fusions.log} for the fusion of conjugacy
classes.

\begin{step-th}
The maximal subgroup $\M_3 \simeq 2^{1+8}.\mathbb{A}_{9}$.
\end{step-th}

Consider the short exact sequence
$1\to2^{1+8}\to2^{1+8}.\mathbb{A}_{9}\to\mathbb{A}_{9}\to 1$. From
Table \ref{tab:altsim}, every non-trivial conjugacy class of
$\mathbb{A}_{9}$ is of type D except the class of the 3-cycles. By
\cite[Lemma 2.7]{AFGV-espo}, every conjugacy class in
$2^{1+8}.\mathbb{A}_{9}$ with representative of order 5, 7, 9, 10,
14, 15, 18, 20, 28, 30, 36 is of type D. Hence the conjugacy
classes \textup{5A, 7A, 10A, 14A, 15A, 15B, 18A, 18B, 20A, 28A,
30A, 30B, 36A, 36B, 36C} of $Th$ are of type D.


\begin{step-th}
The maximal subgroup $\M_2 \simeq 2^{5}.L_5(2)$.
\end{step-th}

 Consider the short exact sequence
$1\to2^{5}\to2^{5}.L_5(2)\to L_5(2)\to 1$. By Lemma
\ref{lem:PSL(5,2)}, in $L_5(2)$ every non-trivial conjugacy class
with representative of order distinct from 2 and 31 is of type D.
Therefore, by \cite[Lemma 2.7]{AFGV-espo}, every conjugacy class
in $\M_2$ with representative of order 3, 6, 12, 21, 24 is of type
D. Hence the classes \textup{3A, 3C, 12D, 21A, 24A, 24B} of $Th$
are of type D.

 On the other hand, each group $L_5(2)$ and
$2^{5}.L_5(2)$ have six classes of elements of order 31; let
$\{\oc_1,\dots, \oc_6 \}$ be these classes in $L_5(2)$. We can
check that for every $i$, $j$, with $1\leq i\neq j\leq 6$, there
exist $r\in \oc_i$ and $s\in \oc_j$ such that $(rs)^2\neq (sr)^2$.
Then the same occurs in the corresponding classes of $
2^{5}.L_5(2)$. Since the fusion of the conjugacy classes from
$2^{5}.L_5(2)$ to $Th$ establish that \textup{31a, 31c, 31d} go to
\textup{31A}, and \textup{31b, 31e, 31f} go to \textup{31B}, then
the conjugacy classes \textup{31A} and \textup{31B} of $Th$ are of
type D. We note that here the classes of this maximal subgroup are
named in lower case letter because they are not necessarily named
as in the \atlas.


\begin{step-th}
The maximal subgroup $\M_{12} \simeq L_2(19).2$.
\end{step-th}

 In this maximal subgroup, the conjugacy classes with
representatives of order $2, 3, 6, 19$ are of type D -- see the file
\texttt{Th/step3.log}. Then, by the fusion of the conjugacy classes,
the conjugacy classes \textup{2A, 3B, 19A} of $Th$ are of type D.


%

\begin{step-th}
The maximal subgroup $\M_5 \simeq (3\times G_{2}(3)):2$.
\end{step-th}

 From Theorem \ref{thm:G2(3)}, the group $G_2(3)$ is of type
D. Thus, the conjugacy class 13A of $Th$ is of type D. On the other hand, by
the Subsubsection \ref{exa:3xG2(3)} and \cite[Lemma 2.8]{AFGV-espo}, 
the conjugacy classes \textup{39A,
39B} of $Th$ are of type D. See the file \texttt{Th/step4.log} for
the fusion of the conjugacy classes $3\times G_2(3)\to (3\times
G_2(3)):2$.



\begin{step-th}
The maximal subgroup $\M_6 \simeq
3.3^{2}.3.(3\times3^{2}).3^{2}:2\mathbb{S}_{4}$.
\end{step-th}

 By the fusion of the conjugacy classes and Algorithm
\ref{alg:tipoD}, we see that the classes \textup{4A, 8A, 12A, 12B,
12C} of $Th$ are of type D. See the file \texttt{Th/step5.log}
for details.


\begin{step-th}
The maximal subgroup $\M_7  \simeq
3^{2}.3^{3}.3^{2}.3^{2}:2\mathbb{S}_{4}$.
\end{step-th}

 By the fusion of the conjugacy classes and Algorithm \ref{alg:tipoD}, we see that
the conjugacy classes \textup{4B, 6A, 6B, 6C, 8B, 9A, 9B, 9C, 24C,
24D, 27A, 27B, 27C} of $Th$ are of type D.  See the file
\texttt{Th/step6.log} for details.



%

\subsection{The Janko group $J_4$}\label{subsect:J4}

See the file \texttt{J4/fusions.log} for the fusion of conjugacy
classes.

\begin{step-j4}
The maximal subgroup $\M_{1}\simeq2^{11}:M_{24}$.
\end{step-j4}

 We use the maximal subgroup $\M_{1}\simeq2^{11}:M_{24}$.
Consider the short exact sequence $1 \to2^{11}\to2^{11}:M_{24}\to
M_{24}\to 1$.  We know that every non-trivial conjugacy class of $M_{24}$ with
representative of order distinct from $23$ is of type D. By  \cite[Lemma 2.7]{AFGV-espo}, every non-trivial conjugacy class in
$2^{11}:M_{24}$ with representative of order distinct from
$2,4,8,16,23$ is of type D. Hence the conjugacy classes
\textup{3A, 5A, 6A, 6B, 6C, 7A, 7B, 10A, 10B, 12A, 12B, 12C, 14A,
14B, 14C, 14D, 15A, 20A, 20B, 21A, 21B, 24A, 24B, 28A, 28B, 30A}
of $J_4$ are of type D. Also, this maximal subgroup has a
primitive permutation representation on $2^{11}$ points. We
construct this primitive group and use Algorithm \ref{alg:tipoD_random} to determine
that the only conjugacy class with representative of order $16$ in $\M_1$ is of type D -- see the file
\texttt{J4/step1.log}. 
Hence, the conjugacy class \textup{16A} of $J_4$ is of type D.



\begin{step-j4}
The maximal subgroup $\M_{4}\simeq2^{3+12}.(\mathbb{S}_5\times L_3(2))$.
\end{step-j4}

 Consider the short exact sequence $1
\to2^{3+12}\to2^{3+12}.(\mathbb{S}_5\times L_3(2))\to (\mathbb{S}_5\times L_3(2))\to
1$. Every non-trivial conjugacy class of $\mathbb{S}_5\times L_3(2)$ with
representative of order distinct from $2,3,4,6,7$ is of type D, see
Subsubsection \ref{exa:S5xL3(2)}. By  \cite[Lemma 2.7]{AFGV-espo}, every
conjugacy class in $2^{3+12}.(\mathbb{S}_5\times L_3(2))$ with
representative of order $35,42$ is of type D. Hence the conjugacy
classes \textup{35A, 35B, 42A, 42B} of $J_4$ are of type D.



\begin{step-j4}
The maximal subgroup $\M_{5}\simeq U_3(11).2$.
\end{step-j4}

 We perform Algorithm \ref{alg:tipoD_random} -- see the
file \texttt{J4/step3.log}. Hence the conjugacy classes
\textup{8A, 8B, 11A, 11B} of $J_4$ are of type D.



\begin{step-j4}
The maximal subgroup $\M_{6}\simeq M_{22}.2$.
\end{step-j4}

 We perform Algorithm \ref{alg:tipoD} -- see the file
\texttt{M22/M22.2.log}. Hence the conjugacy classes \textup{4B,
8C} of $J_4$ are of type D.



\begin{step-j4}
The maximal subgroup $\M_{7}\simeq 11_+^{1+2}:(5\times2\mathbb{S}_4)$.
\end{step-j4}

 We perform Algorithm \ref{alg:tipoD} -- see the file
\texttt{J4/step5.log}. Hence the conjugacy classes \textup{4A,
22A, 22B, 40A, 40B, 44A, 66A, 66B} of $J_4$ are of type D.



\begin{step-j4}
The maximal subgroup $\M_{8}\simeq L_2(32).5$.
\end{step-j4}

 We perform Algorithm \ref{alg:tipoD} -- see the file
\texttt{J4/step6.log}. Hence the conjugacy classes \textup{31A,
31B, 31C, 33A, 33B} of $J_4$ are of type D.


\begin{step-j4}
The maximal subgroup $\M_{9}\simeq L_2(23).2$.
\end{step-j4}

 We perform Algorithm \ref{alg:tipoD} -- see the file
\texttt{J4/step7.log}. Hence the conjugacy classes \textup{2A,
2B, 23A} of $J_4$ are of type D.

\begin{step-j4}
The maximal subgroup $\M_{13}\simeq 37:12$.
\end{step-j4}

 We construct a permutation representation of this
maximal subgroup. By the fusion of the conjugacy classes, the
classes \textup{4a} and \textup{4b} of $\M_{13}$ go to the
conjugacy class \textup{4C} of $J_4$. We find $r$ in \textup{4A},
$s$ in \textup{4B} of $\M_{13}$ such that $(rs)^2 \ne (sr)^2$ --
see the file \texttt{J4/step8.log}. Hence the conjugacy class
\textup{4C} of $J_4$ is of type D. We note that here the classes
of this maximal subgroup are named with lower case letters because
they are not necessarily named as in the \atlas.



\begin{rem}
Not necessarily of type D: \textup{29A, 37A, 37B, 37C, 43A, 43B,
43C}.
\end{rem}

\subsection{The Fischer group $Fi_{23}$}\label{subsect:Fi23}

See the file \texttt{Fi23/fusions.log} for the fusion of conjugacy
classes.

\begin{step-fi23}
The maximal subgroup $\M_{1}\simeq2.Fi_{22}$.
\end{step-fi23}

 Consider the short exact sequence $1\to2\to2.Fi_{22}\to Fi_{22}\to
1$. By Table \ref{tab:alg2}, every non-trivial conjugacy class of $Fi_{22}$
with representative of order distinct from $2,22$ is of type D. By
 \cite[Lemma 2.7]{AFGV-espo}, every conjugacy class in $2.Fi_{22}$
with representative of order 3, 5, 6, 7, 8, 9, 10, 11, 12, 13, 14, 16,
20, 21, 26, 42 is of type D. Then the conjugacy classes \textup{3A, 3B, 3C, 3D, 5A,
6A, 6B, 6C, 6D, 6E, 6F, 6G, 6H, 6I, 6J, 6K, 6L, 6M, 6N, 6O, 7A, 9B, 9C, 9E, 10A, 10B, 10C, 11A,
12A, 12B, 12C, 12E, 12F, 12G, 12H, 12I, 12J, 12K, 12L, 12M, 12N, 12O,
13A, 13B, 14A, 14B, 16A, 16B, 20A, 20B, 21A, 26A, 26B} of $Fi_{23}$ are of type D.

\begin{step-fi23}
The maximal subgroup $\M_4\simeq S_8(2)$.
\end{step-fi23}

 By Lemma \ref{lem:S8(2)_tipoD} and the fusion of the conjugacy classes, the conjugacy classes
\textup{2B, 2C, 4A, 4B, 4C, 4D, 8C, 15A, 15B, 17A} of $Fi_{23}$ are of
type D.

\begin{step-fi23}
The maximal subgroup $\M_5\simeq O_7(3)\times\mathbb{S}_3$.
\end{step-fi23}

 By Lemma \ref{lem:O7(3)_tipoD} and \cite[Lemma 2.8]{AFGV-espo}, every non-trivial conjugacy class of
$O_7(3)\times\mathbb{S}_3$ with representative of order distinct from $2, 3, 6$ is of type D.
Thus by the fusion of the conjugacy classes, the conjugacy classes \textup{9D, 12D, 18A, 18B, 18C, 18E, 18F, 18H, 39A, 39B} of $Fi_{23}$
are of type D.

\begin{step-fi23}
The maximal subgroup $\M_3\simeq 2^2.U_6(2).2$.
\end{step-fi23}

 We perform Algorithm \ref{alg:tipoD_random} in this
maximal subgroup to see that the conjugacy classes \textup{22A,
22B, 22C} of $Fi_{23}$ are of type D. See the file
\texttt{Fi23/step4.log} for details.

\begin{step-fi23}
The maximal subgroup $\M_2\simeq O_8^+(3):\mathbb{S}_3$.
\end{step-fi23}

 We perform Algorithm \ref{alg:tipoD_random} to see that
the conjugacy class \textup{27A} of $Fi_{23}$ is of type D. See
the file \texttt{Fi23/step5.log} for details.

\begin{step-fi23}
The maximal subgroup $\M_{10}\simeq (2^2\times 2^{1+8}).(3\times U_4(2)).2$.
\end{step-fi23}

 We perform Algorithm \ref{alg:tipoD} in this maximal
subgroup to see that the classes \textup{9A, 18D, 18G, 24A, 24B,
24C, 36A} of $Fi_{23}$ are of type D. See
\texttt{Fi23/step6.log} for details.

\begin{step-fi23}
The maximal subgroup $\M_{12}\simeq\mathbb{S}_4\times S_6(2)$.
\end{step-fi23}

 We perform Algorithm \ref{alg:tipoD} in this maximal
subgroup to see that the conjugacy class \textup{36B} $Fi_{23}$ is
of type D. See the file \texttt{Fi23/step7.log} for details.

\begin{step-fi23}
The maximal subgroup $\M_{9}\simeq\mathbb{S}_{12}$.
\end{step-fi23}

 From Table \ref{tab:altsim}, every non-trivial conjugacy
class of $\mathbb{S}_{12}$ with representative of order distinct
of $2, 3, 11$ is of type D. Thus, from the fusion of the conjugacy
classes, the conjugacy classes \textup{8A, 8B, 28A, 30A, 30B, 30C,
35A, 42A, 60A} of $Fi_{23}$ are of type D.

\begin{rem}
Not necessarily of type D: \textup{2A, 23A, 23B}.
\end{rem}

\subsection{The Conway group $Co_1$}
\label{subsect:Co1}

See the file \texttt{Co1/fusions.log} for the fusion of conjugacy
classes.

\begin{step-co1}
The maximal subgroups $\M_{1}\simeq Co_2$ and $\M_4\simeq Co_3$.
\end{step-co1}

 From Table \ref{tab:alg2} and the fusion of the conjugacy classes, the
conjugacy classes \textup{3B, 3C, 4A, 4B, 4C, 4D, 4F, 5B, 5C, 6C, 6D, 6E, 6F,
6G, 6I, 7B, 8B, 8C, 8D, 8E, 9B, 9C, 10D, 10E, 10F, 11A, 14B, 15D, 15E, 16A,
16B, 20B, 20C, 21C, 22A, 28A, 30D, 30E} of $Co_1$ are of type D.

\begin{step-co1}\label{step:Co1max16}
The maximal subgroup $\M_{16}\simeq \mathbb{A}_9\times\mathbb{S}_3$.
\end{step-co1}

 We know that every non-trivial conjugacy class of
$\mathbb{A}_9\times\mathbb{S}_3$ with representative of order distinct from
$2,3$ is of type D, see Subsubsection \ref{exa:A9xS3}. Then, by the fusion of
conjugacy classes, the classes \textup{4E, 5A, 6A, 6B, 6H, 9A, 9B,
10A, 10B, 12L, 12M, 15A, 15C, 30A, 30C} of $Co_1$ are of type D.

\begin{step-co1}
The maximal subgroup $\M_{3}\simeq2^{11}:M_{24}$.
\end{step-co1}

 Consider the short exact sequence
$1\to2^{11}\to2^{11}:M_{24}\to M_{24}\to 1$.  By Table
\ref{tab:alg2} and  \cite[Lemma 2.7]{AFGV-espo}, every non-trivial conjugacy class of $2^{11}:M_{24}$ with
representative of order distinct from $23$ is of type D. By the fusion of the conjugacy classes, the conjugacy class
\textup{3D} of $Co_1$ is of type D.

\begin{step-co1}
The maximal subgroup $\M_{7}\simeq (\mathbb{A}_4\times G_2(4)):2$.
\end{step-co1}

 Note that $\mathbb{A}_4\times G_2(4)$ is
a subgroup of $\M_7$. Now by Subsubsection \ref{exa:A4xG2(4)} and the fusion of
conjugacy classes, the conjugacy classes \textup{14A, 26A} of $Co_1$ are of
type D.

\begin{step-co1}
The maximal subgroup $\M_{14}\simeq(\mathbb{A}_6\times U_3(3)):2$.
\end{step-co1}

 The group $\mathbb{A}_6\times U_3(3)$ is a subgroup of
$Co_1$. It has two conjugacy classes with representatives of order
$28$ and four conjugacy classes with representatives of order
$35$, which are all of type D -- see Subsubsection
\ref{exa:A6xU3(3)}. The maximal subgroup $\M_{14}$ has only one
conjugacy class with representative of order $28$, hence this
conjugacy class is of type D. The fussion of conjugacy classes
says that the conjugacy class \textup{28A} of $\M_{14}$ goes to
\textup{28B} of $Co_1$; thus the conjugacy class \textup{28B} of
$Co_1$ is of type D. On the other hand, the conjugacy class
\textup{35A} of $Co_1$ is of type D since $Co_1$ has only one
conjugacy class with representative of order $35$.

\begin{step-co1}
The maximal subgroup $\M_2\simeq 3.Suz.2$.
\end{step-co1}

 We perform Algorithm \ref{alg:tipoD_random} in this
maximal subgroup to see that the conjugacy classes \textup{7A, 8A,
8F, 10C, 13A, 15B, 21A, 21B, 30B, 33A, 39A, 39B, 42A} of $Co_1$
are of type D. See the file \texttt{Suz/3.Suz.2.log} for the
computations.

\begin{step-co1}
The maximal subgroup $\M_5\simeq 2^{1+8}.O_8^+(2)$.
\end{step-co1}

 By Lemma \ref{lem:O8+(2)_tipoD} and  \cite[Lemma 2.7]{AFGV-espo} the conjugacy classes of $2^{1+8}.O_8^+(2)$ with
representative of order 5, 7, 9, 10, 14, 15, 18, 20, 28, 30, 36, 40, 60 are of type D.
Thus, by the fusion of the conjugacy classes, the conjugacy classes \textup{20A, 36A, 40A, 60A} of $Co_1$ are of type D.

 On the other hand, to study other conjugacy classes of $Co_1$ we use
a script that performs an algorithm similar to Algorithm \ref{alg:tipoD_random}
that we explain briefly here -- see the file \texttt{Co1/step7.log}.  First we
compute all conjugacy classes of $Co_1$ and $\M_5$.  Then we study the
conjugacy classes of $\M_5$ with representative of order 12, 18, 24 and 
discard the corresponding conjugacy classes in $Co_1$ when those are of type D.
At the end of the log file we see that only two conjugacy classes of $Co_1$,
both with representatives of order 12, were not discarded. One of these classes has
centralizer of order 48, the other 72. These are the conjugacy classes
\textup{12L} and \textup{12M} of $Co_1$, and these classes were considered in the Step
\ref{step:Co1max16}.  Therefore, besides the conjugacy classes of the previous
paragraph, the conjugacy classes \textup{12A, 12B, 12C, 12D, 12E, 
12F, 12G, 12H, 12I, 12J, 12K, 18A, 18B, 18C, 24A, 24B, 24C, 24D, 24E, 24F} of 
$Co_1$ are of type D. 


\begin{rem}
Not necessarily of type D: \textup{3A, 23A, 23B}.
\end{rem}

\subsection{The Harada-Norton group $HN$}
\label{subsect:HN}
See the file \texttt{HN/fusions.log} for the fusion of conjugacy
classes.

%

\label{ss:HN}
\begin{step-hn}
The maximal subgroup $\M_{1}\simeq \mathbb{A}_{12}$.
\end{step-hn}

 By Table \ref{tab:altsim}, in this maximal subgroup
every non-trivial conjugacy class with representative of order
distinct from $3,11$ is of type D. Therefore, the conjugacy classes
\textup{2A, 2B, 5A, 5E, 6A, 6B, 6C, 7A, 9A, 15A, 20C, 21A, 30A,
35A, 35B} of $HN$ are of type D.  It remains to prove that the
conjugacy class \textup{11A} of $HN$ is of type D. For that
purpose, let $r=(1\;2\;3\;4\;5\;6\;7\;8\;9\;10\;11)$ and
$s=(1\;2\;3\;4\;5\;6\;7\;8\;9\;11\;10)$ be elements in
$\Alt_{12}$. It is easy to see that $(rs)^2\ne(sr)^2$ and that $r$
and $s$ belong to different conjugacy classes in the group
$\langle r,s\rangle\simeq\Alt_{11}$.  Then, the conjugacy class
\textup{11A} of $HN$ is of type D.

\begin{step-hn}
The maximal subgroup $\M_{11}\simeq M_{12}.2$.
\end{step-hn}

 We perform Algorithm \ref{alg:tipoD} and obtain that
every non-trivial conjugacy class is of type D -- see the file
\texttt{M12/M12.2.log}. By the fusion of the conjugacy
classes, the conjugacy classes \textup{3A, 3B, 4A, 4B, 4C, 12C} of
$HN$ are of type D.

\begin{step-hn}
The maximal subgroup $\M_2\simeq 2.HS.2$.
\end{step-hn}

 We perform Algorithm \ref{alg:tipoD_random} and obtain
that every conjugacy class with representative of order $5, 8, 10,
12, 14, 20, 22, 40$ is of type D -- see the file
\texttt{HS/2.HS.2.log}. Therefore, the conjugacy classes
\textup{5B, 8A, 8B, 10A, 10B, 10C, 10F, 10G, 10H, 12A, 12B, 14A,
20A, 20B, 22A, 40A, 40B} of $HN$ are of type D.

\begin{step-hn}
The maximal subgroup $\M_{14}\simeq 3^{(1+4)}:4.\mathbb{A}_5$.
\end{step-hn}

 We use Algorithm \ref{alg:tipoD} and obtain that every
conjugacy class with representative of order $5, 10, 15, 20, 30$
is of type D -- see the file \texttt{HN/step4.log}.  Therefore,
the conjugacy classes \textup{5C, 5D, 10D, 10E, 15B, 15C, 20D,
20E, 30B, 30C} of $HN$ are of type D.

\begin{step-hn}
The maximal subgroup $\M_{3}\simeq U_3(8).3_1$.
\end{step-hn}

 We perform Algorithm \ref{alg:tipoD_random} and obtain
that every conjugacy class with representative of order $19$ is of
type D -- see the file \texttt{HN/step5.log}. Therefore, the
conjugacy classes \textup{19A, 19B} of $HN$ are of type D.

\begin{step-hn}
The maximal subgroup $\M_{10}\simeq 5^{2+1+2}.4.\Alt_5$.
\end{step-hn}

We perform Algorithm \ref{alg:tipoD} and obtain that
every conjugacy class with representative of order $25$ is of type
D -- see the file \texttt{HN/step6.log}. Therefore, the
conjugacy classes \textup{25A, 25B}  of $HN$ are of type D.

\subsection{The Fischer group $Fi_{24}'$}
\label{subsect:Fi24'}

The log files concerning this group are stored in the folder 
\texttt{F3+}. See the file \texttt{F3+/fusions.log} for the fusion of
the conjugacy classes. The list of (representatives of conjugacy
classes of) maximal subgroups of $Fi_{24}'$ can be found in
\cite{LW91}.


\begin{step-fi24'}
The maximal subgroup $\M_{1}\simeq Fi_{23}$.
\end{step-fi24'}
 By Section \ref{subsect:Fi23}, we know that every non-trivial conjugacy class
of $Fi_{23}$ with representative of order distinct from $2,23$ is of type D.
Hence, the conjugacy classes \textup{3A, 3B, 3C, 3D, 4A, 4B, 4C,
5A, 6A, 6B, 6C, 6D, 6E, 6F, 6G, 6H, 6I, 6J, 7A, 8A, 8B, 9A, 9B,
9C, 9E, 9F, 10A, 10B, 11A, 12A, 12B, 12C, 12D, 12E, 12F, 12G, 12H,
12K, 12L, 12M, 13A, 14A, 15A, 15C, 16A, 17A, 18A, 18B, 18C, 18D,
18E, 18F, 20A, 21A, 22A, 24A, 24B, 24E, 26A, 27A, 28A, 30A, 30B,
35A, 36C, 36D, 39A, 39B, 42A, 60A} of $Fi_{24}'$ are of type D.

\begin{step-fi24'}
The maximal subgroups $\M_{13}\simeq He.2$ and $\M_{14}\simeq He.2$.
\end{step-fi24'}

  We perform Algorithm \ref{alg:tipoD_maximales} and
obtain that every non-trivial conjugacy class is of type D -- see
the file \texttt{He/He.2.log}. Then the conjugacy classes
\textup{2A, 2B, 3E, 6K, 7B, 12I, 12J, 14B, 21B, 21C, 21D, 24C,
24D, 42B, 42C} of $Fi_{24}'$ are of type D.

\begin{step-fi24'}
The maximal subgroup $\M_{4}\simeq O_{10}^{-}(2)$.
\end{step-fi24'}

 By Lemma \ref{lem:O10-(2)_tipoD} and the fusion of the conjugacy classes, the conjugacy classes
\textup{8C, 15B, 18G, 18H, 20B} of $Fi_{24}'$ are of type D.

\begin{step-fi24'}
The maximal subgroup $\M_{5}\simeq 3^7.O_7(3)$.
\end{step-fi24'}

 We consider the
short exact sequence $1\to3^7\to3^7.O_7(3)\to O_7(3)\to 1$. By
Lemma \ref{lem:O7(3)_tipoD}, every non-trivial conjugacy class of $O_7(3)$
with representative of order distinct from $2$ is of type D.  \cite[Lemma 2.7]{AFGV-espo}, every non-trivial conjugacy class of $\M_5$
with representative of order distinct from $24, 36, 45$ is of type D.
Therefore, the conjugacy classes \textup{24F, 24G,
36A, 36B, 45A, 45B} of $Fi_{24}'$ are of type D.

\begin{step-fi24'}
The maximal subgroup $\M_{20}\simeq\Alt_6\times L_2(8):3$.
\end{step-fi24'}

 We perform Algorithm \ref{alg:tipoD} and obtain that
every conjugacy class with representative of order $9$ is of type
D -- see the file \texttt{F3+/step5.log}. Therefore, the
conjugacy class \textup{9D}  of $HN$ is of type D.

\begin{rem}
Not necessarily of type D: \textup{23A, 23B, 27B, 27C, 29A, 29B, 33A, 33B, 39C, 39D}.
\end{rem}

\subsection{The Baby Monster group $B$}\label{subsect:B}
For this group we compute the fusion of the conjugacy classes in each step.

\label{ss:B1}
\begin{step-b}
The maximal subgroup $\M_{2}\simeq 2^{1+22}.Co_2$.
\end{step-b}

 We consider the short exact sequence $1\to2^{1+22}\to2^{1+22}.Co_2\to Co_2\to
 1$. By Table \ref{tab:alg2}, every non-trivial conjugacy class of $Co_2$ with
 representative of order distinct from $2,23$ is of type D.  Therefore, by
 \cite[Lemma 2.7]{AFGV-espo}, the conjugacy classes \textup{5A, 5B, 6A, 6B, 6C,
 6D, 6E, 6F, 6G, 6H, 6I, 6J, 6K, 7A, 9A, 9B, 10A, 10B, 10C, 10D, 10E, 10F, 11A,
 12A, 12B, 12C, 12D, 12E, 12F 12G, 12H, 12I, 12J, 12K, 12L, 12M, 12N, 12O, 12P,
 12Q, 12R, 12S, 14A, 14B, 14C, 14D, 14E, 15A, 15B, 18F, 20A, 20B, 20C, 20D,
 20E, 20F, 20G, 20H, 20I, 20J, 24A, 24B, 24C, 24D, 24E, 24F, 24G, 24H, 24I,
 24J, 24K, 24M, 28A, 28B, 28C, 28D, 28E, 30A, 30B, 30C, 30D, 30E, 30F, 30G,
 30H, 36A, 40A, 40B, 40C, 40D, 44A, 48A, 48B, 56A, 56B} of $B$ are of type D.
 See the file \texttt{B/step1.log} for the fusion of the conjugacy classes
 $2^{1+22}.Co_2\to B$.

\label{ss:B2}
\begin{step-b}
The maximal subgroup $\M_{3}\simeq Fi_{23}$.
\end{step-b}

 By Subsection \ref{subsect:Fi23}, every non-trivial conjugacy class of $Fi_{23}$ distinct
of \textup{2A, 23A, 23B} is of type D.  Therefore, the conjugacy classes
\textup{2B, 2D, 3A, 3B, 4D, 4E, 4G, 4H, 8J, 8K, 13A, 16G, 17A, 18A, 18B, 18C,
18D, 18E, 21A, 22A, 22B, 24L, 26B, 27A, 35A, 36B, 36C, 39A, 42B} of $B$
are of type D. See the file \texttt{B/step2.log} for the fusion of conjugacy
classes $Fi_{23}\to B$.

\begin{step-b}
The maximal subgroup $\M_{18}\simeq \s_5 \times M_{22}:2$.
\end{step-b}

 We perform Algorithm \ref{alg:tipoD} and obtain that
every conjugacy class of $\M_{18}$ with representative of order
$4, 8, 33, 42, 55, 66, 70$ is of type D. Hence the conjugacy
classes \textup{4A, 4B, 4C, 4F, 8B, 8C, 8E, 8I, 33A, 42C, 55A,
66A, 70A} of $B$ are of type D. See the file \texttt{B/step3.log}
for the fusion of conjugacy classes $\Sym_5\times M_{22}:2\to B$
and the computations.

\begin{step-b}
The maximal subgroup $\M_{4}\simeq 2^{9+16}.S_8(2)$.
\end{step-b}

 By Lemma \ref{lem:S8(2)_tipoD}, every non-trivial conjugacy class, except
2A, 2B, 3A is of type D. By  \cite[Lemma 2.7]{AFGV-espo} and the fusion of conjugacy
classes -- see the file \texttt{B/step4.log} --  the conjugacy
classes \textup{34B, 34C, 40E, 42A, 60A, 60B, 60C} of $B$ are of type D.

\begin{step-b}
The maximal subgroup $\M_{16}\simeq \mathbb{S}_4\times{}^2F_4(2)$.
\end{step-b}

 See the file \texttt{B/step5.log} for the fusion of the conjugacy classes from this
maximal subgroup into $B$. By Lemma \ref{lem:2F4(2)_tipoD} and \cite[Lemma 2.8]{AFGV-espo}, we deduce that every non-trivial conjugacy class of $\M_{16}$
with representative of order distinct from $2, 3, 4$ is of type D. Thus, the conjugacy classes \textup{4I, 8A, 8D, 8F, 8G,
8H, 8L, 8N, 12T, 16A, 16B, 16E, 16F, 26A, 52A} of $B$ are of type D.

\begin{step-b}
The maximal subgroup $\M_{28}\simeq L_2(17).2$.
\end{step-b}

 We perform Algorithm \ref{alg:tipoD} and obtain that
every conjugacy class of $\M_{28}$ with representative of order
$8, 16$ is of type D.  Then, by the fusion of the conjugacy
classes, the conjugacy classes \textup{8M, 16H} of $B$ are of type
D. See the file \texttt{B/step6.log} for the fusion of the
conjugacy classes and the computations.

\begin{step-b}
The maximal subgroup $\M_{5}\simeq Th$.
\end{step-b}

 By Subsection \ref{subsect:Th}, every non-trivial conjugacy class in $Th$ is of type D.
Therefore, by the fusion of the conjugacy classes -- see the file \texttt{B/step7.log} --
the conjugacy classes \textup{4J, 19A, 24N, 31A, 31B} of $B$ are of type D.

\begin{step-b}
The maximal subgroup $\M_{30}\simeq 47:23$.
\end{step-b}

 We use this maximal subgroup to show that the classes \textup{23A}, \textup{23B} of $B$ are of
type D. Indeed, the classes \textup{23a} and \textup{23d} of
$\M_{30}$ fuse into the class \textup{23A} of $B$ and there exist
elements $r$ and $s$ in the classes \textup{23a} and \textup{23d}
of $\M_{30}$, respectively, such that $(rs)^2\neq (sr)^2$. Then
the class \textup{23A} of $B$ is of type D. On the other hand, the
classes \textup{23u} and \textup{23v} of $\M_{30}$ fuse into the
class \textup{23B} of $B$ and there exist elements $r$ and $s$ in
the classes \textup{23u} and \textup{23v} of $\M_{30}$,
respectively, such that $(rs)^2\neq (sr)^2$. Then the class
\textup{23B} of $B$ is of type D. See the file
\texttt{B/step8.log} for the details. We note that here the
classes of this maximal subgroup are named in lower case letter
because they are not necessarily named as in the \atlas.

\begin{step-b}
The maximal subgroup $\M_1\simeq 2.(^2E_6(2)):2$.
\end{step-b}

 We know that $2\times U_3(8)$ is a subgroup of $B$,
because $2\times U_3(8)$ is a subgroup of $\M_1$, see \cite{W}. In
$U_3(8)$ there exist non-conjugate elements $r$ and $s$, both of
order $19$, such that $(rs)^2\ne(sr)^2$ -- see the file
\texttt{B/step9.log}. Since $B$ has only one conjugacy class of
elements of order $38$, the conjugacy class \textup{38A} of $B$ is
of type D.

\begin{rem}
Not necessarily of type D: \textup{2A, 
16C, 16D, 32A, 32B, 32C, 32D, 34A, 46A, 46B, 47A, 47B}
\end{rem}

\subsection{The Monster group $M$}\label{subsect:M}

For the known maximal subgroups of $M$, we use the order in which
they appear listed in the \atlas.

%

%

\begin{step-m}
The maximal subgroup $\M_{1}\simeq 2.B$.
\end{step-m}

 By Subsection \ref{subsect:B} and  \cite[Lemma
2.7]{AFGV-espo}, the conjugacy classes \textup{3A, 3B, 5A, 5B, 6A,
6B, 6C, 6D, 6E, 7A, 8A, 8B, 8C, 8D, 8E, 8F, 9A, 9B, 10A, 10B, 10C,
10D, 10E, 11A, 12A, 12B, 12C, 12E, 12F, 12G, 12H, 12I, 13A, 14A,
14B, 17A, 18A, 18B, 18C, 18D, 18E, 19A, 21A, 22A, 22B, 23A, 23B,
24A, 24B, 24C, 24D, 24F, 24G, 24H, 24I, 25A, 26A, 27A, 28A, 28B,
28C, 33B, 35A, 36A, 36B, 36C, 36D, 38A, 39A, 40A, 40B, 40C, 40D,
42A, 44A, 44B, 48A, 50A, 52A, 54A, 55A, 66A, 70A, 78A, 84A, 104A,
104B, 110A} of $M$ are of type D. See the file
\texttt{M/step1.log} for the fusion of conjugacy classes $2.B\to
M$.

\begin{step-m}
The maximal subgroup $\M_{2}\simeq 2^{1+24}.Co_1$.
\end{step-m}

 By Subsection \ref{subsect:Co1} and
 \cite[Lemma 2.7]{AFGV-espo}, the conjugacy classes \textup{7B, 13B, 14C, 15A, 15B, 15C,
15D, 20A, 20B, 20C, 20D, 20E, 20F, 21B, 21D, 26B, 28D, 30A, 30B,
30C, 30D, 30E, 30F, 30G, 33A, 35B, 39C, 39D, 42B, 42C, 42D, 52B,
56A, 56B, 56C, 60A, 60B, 60C, 60D, 60E, 60F, 66B, 70B, 78B, 78C,
84B, 88A, 88B} of $M$ are of type D. See the file
\texttt{M/step2.log} for the fusion of conjugacy classes
$2^{1+24}.Co_1\to M$.

\begin{step-m}
The maximal subgroup $\M_{9}\simeq \mathbb{S}_3\times Th$.
\end{step-m}

 By Subsection \ref{subsect:Th} and \cite[Lemma
2.8]{AFGV-espo}, every non-trivial conjugacy class in this maximal
subgroup with representative of order distinct from $2, 3$ is of
type D. Moreover, in this maximal subgroup, there is only one
conjugacy class with representative of order 3 that is not of type
D: that corresponding to the 3-cycles in $\Sym_3$. Hence, the
conjugacy classes \textup{3C, 4A, 4D, 6F, 12D, 12J, 21C, 24E, 24J,
27B, 31A, 31B, 39B, 57A, 62A, 62B, 84C, 93A, 93B} of $M$ are of
type D. See the file \texttt{M/step3.log} for the fusion of
conjugacy classes $\s_3\times Th\to M$

\begin{step-m}
The maximal subgroup $\M_{40}\simeq L_2(29).2\simeq
\mathbf{PGL}(2,29)$.
\end{step-m}

 We use a representation of this maximal subgroup inside
$\Sym_{30}$ given in the \textsf{ATLAS}, see \cite{BW}.  This maximal
subgroup has only one conjugacy class of elements of order 29
which is of type D. Therefore, the conjugacy class \textup{29A} of
$M$ is of type D.  See the file \texttt{M/step4.log} for details.

\begin{step-m}
The maximal subgroup $\M_{23}\simeq (L_3(2)\times S_4(4):2).2$.
\end{step-m}

 We use a representation of this maximal subgroup inside
$\Sym_{184}$ given in the \textsf{ATLAS}, see \cite{BW}.  In this
maximal subgroup the conjugacy classes with representatives of
order 16, 34, 51, 68, 119 are of type D. Therefore, the conjugacy
classes 16A, 16B, 16C, 34A, 51A, 68A, 119A, 119B of $M$ are of
type D. See the file \texttt{M/step5.log} for details.

\begin{step-m}
The maximal subgroup $\M_{21}\simeq (\Alt_5\times U_3(8):3_1):2$.
\end{step-m}

 We use a representation of this maximal subgroup inside
$\Sym_{518}$ given in the \textsf{ATLAS}, see \cite{BW}. In this
maximal subgroup, the conjugacy classes with representatives of
order 95 are of type D. Therefore, the conjugacy classes
\textup{95A, 95B} of $M$ are of type D.  See the file
\texttt{M/step6.log} for details.

\begin{step-m}
The maximal subgroup $\M_{3}\simeq 3.Fi_{24}$.
\end{step-m}

 Let $H=Fi_{24}$. From the \atlas~we know that the group
$K=\Sym_5\times\Sym_9$ is a maximal subgroup of $H$. From
\ref{exa:S5xS9} every conjugacy class of $K$ with representative
of order 4 is of type D. Also, from the fusion of conjugacy
classes $K\to H$, every conjugacy class of $Fi_{24}$ with
representative of order 4 is of type D.  Therefore, by \cite[Lemma
2.7]{AFGV-espo}, the conjugacy classes \textup{4A, 4B, 4C, 4D} of
$M$ are of type D. Also, since $Fi_{24}'$ is a maximal subgroup of
$Fi_{24}$, the conjugacy classes of $Fi_{24}$ with representatives
of order 15, 45 are of type D -- see the Subsection
\ref{subsect:Fi24'}. Therefore, by  \cite[Lemma 2.7]{AFGV-espo}
and the fusion of conjugacy classes $Fi_{24}'\to Fi_{24}$, the
conjugacy class \textup{45A} of $M$ is of type D. For details
about these observations and the fusion of conjugacy classes
$3.Fi_{24}\to M$ see the file \texttt{M/step7.log}.

\begin{step-m}
The subgroup $HN$.
\end{step-m}

 Since $(\mathbb{D}_{10}\times HN).2$ is a maximal
subgroup of $M$, $HN$ is a subgroup of $M$. Therefore, by
Subsection \ref{subsect:HN} and the fusion of the conjugacy
classes $HN\to M$, the conjugacy clases \textup{2A, 2B} of $M$ are
of type D. For the fusion of the conjugacy classes see the file
\texttt{M/step8.log}.

\begin{step-m}
The subgroup $2\times47:23$
\end{step-m}

We use this subgroup to show that the classes \textup{46C}, 
\textup{46D} of $M$ are of
type D. Indeed, the classes \textup{46A} and \textup{46B} of
$2\times47:23$ fuse into the class \textup{46C} of $M$ and there exist
elements $r$ and $s$ in the classes \textup{46A} and \textup{46B}
of $2\times47:23$, respectively, such that $(rs)^2\neq (sr)^2$. Then
the class \textup{46C} of $M$ is of type D. The conjugacy class
\textup{46D} is treated analogously. 
See the file \texttt{M/step9.log} for the details. 

\begin{rem}
Not necessarily of type D: \textup{32A, 32B, 41A, 46A, 46B, 
47A, 47B, 59A, 59B, 69A, 69B, 71A, 71B, 87A, 87B, 92A, 92B,
94A, 94B}.
\end{rem}

\section*{Appendix. Real and quasi-real conjugacy classes}
\label{sect:appendix}

In this appendix we list all real and quasi-real conjugacy classes of the
groups studied in \cite{AFGV-espo}. The information about real conjugacy
classes of a given group $G$ is easy to obtain from the character table of $G$
using the \gap~function \texttt{RealClasses}. Similarly, with the \gap~function
\texttt{PowerMaps} it is easy to determine the quasi-real conjugacy classes of
a given group.

The function \texttt{QuasiRealConjugacyClasses} returns the list of quasi-real
conjugacy classes and its type.

\begin{verbatim}
     gap> QuasiRealClasses := function( ct )
     >  local nc, oc, a, b, p, c, j, rc;
     >
     >  nc := NrConjugacyClasses(ct);
     >  oc := OrdersClassRepresentatives(ct);
     >  rc := RealClasses(ct);
     >
     >  a := [];
     >  b := [];
     >
     >  for c in [1..nc] do
     >    if not c in rc then
     >      for j in [2..oc[c]-2] do
     >        p := PowerMap(ct, j);
     >        if p[c] = c then
     >          if j-1 mod oc[c] <> 0 then
     >            if not c in b then
     >              Add(a, [c, j]);
     >              Add(b, c);
     >            fi;
     >          fi;
     >        fi;
     >      od;
     >    fi;
     >  od;
     >  return a;
     >end;
\end{verbatim}

\subsubsection*{The group $L_5(2)$}

The conjugacy classes \textup{7A, 7B, 15A, 15B, 21A, 21B, 31A,
31B, 31C, 31D, 31E, 31F} are quasi-real of type $j=2$, and the
classes \textup{14A, 14B} are quasi-real of type $j=9$. The
remaining conjugacy classes are real.


\subsubsection*{The groups $O^+_8(2)$, $S_6(2)$ and $S_8(2)$}

In these groups every conjugacy class is real.

\subsubsection*{The group $O^-_{10}(2)$}

The conjugacy classes \textup{3B, 3C, 6B, 6C, 6F, 6G, 6H, 6I, 6L,
6M, 6T, 6U} are neither real nor quasi-real. The conjugacy classes
\textup{11A, 11B, 35A, 35B} are quasi-real of type $j=3$, the
classes \textup{9B, 9C, 15B, 15C, 15F, 15G, 33A, 33B, 33C, 33D}
are quasi-real of type $j=4$, the classes \textup{12B, 12C, 12E,
12F, 12I, 12J, 12N, 12O, 12R, 12S, 18A, 18B, 18C, 18D, 24C, 24D,
30B, 30C} are quasi-real of type $j=7$. The remaining conjugacy
classes are real.



\subsubsection*{The group $G_2(4)$}

The conjugacy classes \textup{12B, 12C} are quasi-real of type
$j=7$.  The remaining conjugacy classes are real.


\subsubsection*{The Tits group}

The conjugacy classes \textup{8A, 8B} are quasi-real of type
$j=5$, the conjugacy classes \textup{16A, 16B, 16C, 16D} are
quasi-real of type $j=9$. The remaining conjugacy classes are
real.

\subsubsection*{The Mathieu groups}

In any of the Mathieu simple groups, every conjugacy class is real or
quasi-real. See Table \ref{tab:Mathieu_classes} for the details concerning not
real but quasi-real conjugacy classes.

\begin{table}[ht]
\caption{Mathieu groups: quasi-real classes.}
\label{tab:Mathieu_classes}
\begin{center}
\begin{tabular}{|c|c|c|}
\hline
        & Classes               & Type\\\hline
$M_{11}$    & \textup{8A, 8B, 11A, 11B}          & $j=3$\\\hline
$M_{12}$    & \textup{11A, 11B}              & $j=3$\\\hline
$M_{22}$    & \textup{7A, 7B}                & $j=2$\\
        & \textup{11A, 11B}              & $j=3$\\\hline
$M_{23}$    & \textup{7A, 7B, 15A, 15B, 23A, 23B}        & $j=2$\\
        & \textup{11A, 11B}              & $j=3$\\
        & \textup{14A, 14B}              & $j=9$\\\hline
$M_{24}$    & \textup{7A, 7B, 15A, 15B, 21A, 21B, 23A, 23B}  & $j=2$\\
        & \textup{14A, 14B}              & $j=9$\\\hline
\end{tabular}
\end{center}
\end{table}

\subsubsection*{The Conway groups}

In the Conway groups $Co_{1}$, $Co_{2}$ and $Co_{3}$ every conjugacy class is
real or quasi-real. The quasi-real not real conjugacy classes are listed in
Table \ref{tab:Conway_classes}.

\begin{table}[ht]
\caption{Conway groups: quasi-real classes.}
\label{tab:Conway_classes}
\begin{center}
\begin{tabular}{|c|c|c|}
\hline
        & Classes               & Type\\\hline
$Co_{1}$    & \textup{23A, 23B, 39A, 39B}            &
$j=2$\\\hline
$Co_{2}$    & \textup{15B, 15C, 23A, 23B}            & $j=2$\\
        & \textup{14B, 14C}              & $j=9$\\
        & \textup{30B, 30C}              & $j=17$\\\hline
$Co_{3}$    & \textup{23A, 23B}              & $j=2$\\
        & \textup{11A, 11B, 20A, 20B, 22A, 22B}      & $j=3$\\\hline
\end{tabular}
\end{center}
\end{table}

\subsubsection*{The Janko groups}
In the Janko groups $J_1$ and $J_2$ every conjugacy class is real. In the Janko
group $J_3$ the conjugacy classes 19A, 19B are quasi-real of type $j=4$ and the
remaining conjugacy classes are real. In the group $J_4$ every conjugacy class
is real, with the exceptions of the following classes which are quasi-real:
\begin{enumerate}
\item \textup{7A, 7B, 21A, 21B, 35A, 35B} (of type $j=2$);
\item \textup{14A, 14B, 14C, 14D, 28A, 28B} (of type $j=9$);
\item \textup{42A, 42B} (of type $j=11$).
\end{enumerate}

\subsubsection*{The Fischer groups}

In the Fischer groups $Fi_{22}$, $Fi_{23}$ and $Fi_{24}'$ every conjugacy class
is real or quasi-real. The quasi-real not real conjugacy classes are listed in
Table \ref{tab:Fischer_classes}.

\begin{table}[ht]
\caption{Fischer groups: quasi-real classes.}
\label{tab:Fischer_classes}
\begin{center}
\begin{tabular}{|c|c|c|}
\hline
        & Classes               & Type\\\hline
$Fi_{22}$   & \textup{11A, 11B, 16A, 16B, 22A, 22B}      & $j=3$\\
            & \textup{18A, 18B}                          & $j=7$\\\hline
$Fi_{23}$   & \textup{16A, 16B, 22B, 22C}            & $j=3$\\
        & \textup{23A, 23B}              & $j=2$\\\hline
$Fi_{24}'$  & \textup{23A, 23B}              & $j=2$\\
        & \textup{18G, 18H}                          & $j=7$\\\hline
\end{tabular}
\end{center}
\end{table}

\subsubsection*{The Highman-Sims group}

The conjugacy classes 11A, 11B, 20A, 20B are quasi-real of type
$j=3$. The remaining conjugacy classes are real.

\subsubsection*{The Lyons group}

The conjugacy classes \textup{11A, 11B, 22A, 22B} are quasi-real
of type $j=3$, the conjugacy classes \textup{33A, 33B} are
quasi-real of type $j=4$. The remaining conjugacy classes are
real.

\subsubsection*{The Harada-Norton group}

Every conjugacy class is real, with the exceptions of the
following classes which are quasi-real:
\begin{enumerate}
\item \textup{19A, 19B} (of type $j=4$);
\item \textup{35A, 35B} (of type $j=3$);
\item \textup{40A, 40B} (of type $j=7$).
\end{enumerate}

\subsubsection*{The Held group}

Every conjugacy class is real, with the exceptions of the
following classes which are quasi-real:
\begin{enumerate}
\item \textup{7A, 7B, 7D, 7E, 21C, 21D} (of type $j=2$);
\item \textup{14A, 14B, 14C, 14D, 28A, 28B} (of type $j=9$).
\end{enumerate}

\subsubsection*{The MacLaughlin group}

Every conjugacy class is real, with the exceptions of the
following classes which are quasi-real:

\begin{enumerate}
\item \textup{7A, 7B, 15A, 15B} (of type $j=2$);
\item \textup{11A, 11B} (of type $j=3$);
\item \textup{9A, 9B} (of type $j=4$);
\item \textup{14A, 14B} (of type $j=9$);
\item \textup{30A, 30B} (of type $j=17$).
\end{enumerate}

\subsubsection*{The O'Nan group}

Every conjugacy class is real, with the exceptions of the
following classes which are quasi-real:
\begin{enumerate}
\item \textup{31A, 31B} (of type $j=2$);
\item \textup{20A, 20B} (of type $j=3$).
\end{enumerate}

\subsubsection*{The Rudvalis group $Ru$}

The conjugacy classes \textup{16A, 16B} are quasi-real of type
$j=5$. The remaining conjugacy classes are real.

\subsubsection*{The Suzuki group $Suz$}

The conjugacy classes \textup{6B, 6C}, with centralizers of size
$1296$ are neither real nor quasi-real. The classes \textup{9A,
9B} are quasi-real of type $j=4$, and the classes \textup{18A,
18B} are quasi-real of type $j=7$. The remaining conjugacy classes
are real.

\subsubsection*{The Thompson group}

Every conjugacy class is real, with the exceptions of the
following classes which are quasi-real:
\begin{enumerate}
\item \textup{15A, 15B, 31A, 31B, 39A, 39B} (of type $j=2$);
\item \textup{27B, 27C} (of type $j=4$);
\item \textup{24C, 24D} (of type $j=5$);
\item \textup{12A, 12B, 24A, 24B, 36B, 36C} (of type $j=7$);
\item \textup{30A, 30B} (of type $j=17$).
\end{enumerate}

\subsubsection*{The Baby Monster group}

The conjugacy classes \textup{23A, 23B, 31A, 31B, 47A, 47B} are
quasi-real of type $j=2$, the classes \textup{30G, 30H} are
quasi-real of type $j=17$, the classes \textup{32C, 32D, 46A, 46B}
are quasi-real of type $j=3$. The remaining conjugacy classes are
real.

\subsubsection*{The Monster group}

The conjugacy classes \textup{23A, 23B, 31A, 31B, 39C, 39D, 47A,
47B, 69A, 69B, 71A, 71B, 87A, 87B, 93A, 93B, 95A, 95B, 119A, 119B}
are quasi-real of type $j=2$. The classes \textup{40C, 40D, 44A,
44B, 46A, 46B, 46C, 46D, 56B, 56C, 59A, 59B, 88A, 88B, 92A, 92B,
94A, 94B, 104A, 104B} are quasi-real of type $j=3$.  The classes
\textup{62A, 62B, 78B, 78C} are quasi-real of type $j=5$. The
remaining conjugacy classes are real.

\begin{acknowledgements*}
We are very grateful to Alexander Hulpke, John Bray, Robert Wilson
and very specially to Thomas Breuer for answering our endless
questions on \textsf{GAP}. We also thank Enrique Tobis and the people from \textsf{shiva} and \textsf{ganesh}
who allowed us to use their computers.
\end{acknowledgements*}

\end{document}